\documentclass[12pt]{amsart}

\usepackage{amsmath,amsthm,amsfonts,verbatim,dsfont}
\usepackage{xypic}
\xyoption{all} \xyoption{matrix}

\def\wt{\widetilde}
\def\ra{\overline}
\def\m{\mathfrak{m}}
\def\Z{\mathcal{Z}}

\def\Aut{\mathrm{Aut}}
\def\Hom{\mathrm{Hom}}

\def\id{\mathrm{id}}

\def\op{\mathrm{op}}

\DeclareMathOperator{\even}{even}

\def\x{\mathbf{x}}
\def\i{\mathbf{i}}

\def\a{\mathbf{a}}

\def\v{\mathbf{v}}

\DeclareMathOperator{\cent}{\mathcal{Z}} \DeclareMathOperator{\Ann}{\mathsf{Ann}}
\def\I{\mathds{I}}
\def\J{\mathds{J}}

\def\C{\mathbb{C}}

\def\HH{\mathrm{HH}}

\def\H{\mathrm{H}}

\def\ot{\otimes}
\def\cl#1{{\langle #1\rangle}}

\theoremstyle{plain}
\newtheorem{teo}{Theorem}[section]
\newtheorem{lem}[teo]{Lemma}
\newtheorem{coro}[teo]{Corollary}
\newtheorem{prop}[teo]{Proposition}
\theoremstyle{definition}

\newtheorem{rem}[teo]{Remark}
\theoremstyle{definition}
\newtheorem{defi}[teo]{Definition}

\title[The cohomology of monogenic extensions]{The cohomology of monogenic extensions in
the noncommutative setting}

\begin{document}

\author{Marco Farinati}
\address{Departamento de Matem\'atica\\ Facultad de Ciencias Exactas y Naturales, Pabell\'on
1 - Ciudad Universitaria\\ (1428) Buenos Aires, Argentina.} \curraddr{} \email{mfarinat@dm.uba.ar}
\thanks{Supported by UBACYT X169 and CONICET: PIP 5099}

\author{Jorge A. Guccione}
\address{Departamento de Matem\'atica\\ Facultad de Ciencias Exactas y Naturales, Pabell\'on
1 - Ciudad Universitaria\\ (1428) Buenos Aires, Argentina.} \curraddr{} \email{vander@dm.uba.ar}
\thanks{Supported by PICT 12330, UBACYT X0294 and CONICET}

\author{Juan J. Guccione}
\address{Departamento de Matem\'atica\\ Facultad de Ciencias Exactas y Naturales\\
Pabell\'on 1 - Ciudad Universitaria\\ (1428) Buenos Aires, Argentina.} \curraddr{}
\email{jjgucci@dm.uba.ar}
\thanks{Supported by PICT 12330, UBACYT X0294 and CONICET}

\date{}


\dedicatory{}


\begin{abstract} We extend the notion of monogenic extension to the noncommutative setting, and
we study the Hochschild cohomology ring of such an extension. As an aplication we complete the
computation of the cohomology ring of the rank one Hopf algebras begun in \cite{B-W}.
\end{abstract}

\maketitle

\tableofcontents
\section*{Introduction}
Let $k$ be a field, $G$ a finite group whose order is relative prime to the characteristic of $k$
and $C = k[x]/\langle x^n\rangle$, where $n\ge 2$. Let $\chi\colon G\to k^{\times}$ be a
character. The group $G$ acts on $C$ via ${}^gx = \chi(g) x$ for all $g\in G$. Consider the
corresponding skew algebra $H = C\# k[G]$. As a vector space $H = C\ot k[G]$ and the
multiplication is
$$
(a\ot g)(b\ot h) = a{}^gb\ot gh\quad\text{for all $a,b\in C$ and $g,h\in G$.}
$$
Assume there is a central element $g_1\in G$ such that $\chi(g_1)$ is a primitive $n$-th root of
$1$. Then $H$ is a Hopf algebra with coproduct $\Delta$, counit $\epsilon$ and antipode $S$
defined by
\begin{alignat*}{2}
&\Delta(x) = x\ot 1 + g_1\ot x, &&\qquad \Delta(g)=g\ot g,\\
&\epsilon(x)=0,&&\qquad\epsilon(g)=1,\\
& S(x)= -g_1^{-1}x,&& \qquad S(g)=g^{-1},
\end{alignat*}
for all $g \in G$. These algebras are examples of Rank one Hopf algebras. A generalization of Taft
algebras defined in \cite{A-S} for abelian groups $G$ and generalized further in \cite{K-R} (but
with the opposite coproduct). The cohomology ring of $H$ was computed in \cite{B-W}, where was
also possed the problem of compute the cohomology ring of the other rank one Hopf algebras. That
is, those ones in which the relation $x^n=0$ is replaced by $x^n=g_1^n-1$ (see \cite{K-R}).

\smallskip

The aim of this paper is to solve this problem. We carry out this task by computing the cohomology
ring of the monogenic extensions $K[x,\alpha]/\langle f\rangle$ of a separable $k$-algebra $K$
(here $K[x,\alpha]$ is an endomorphism type Ore extension and $f\in K[x,\alpha]$ is a monic
polynomial satisfying suitable hypothesis) and noting that each rank one Hopf algebras is such an
extension (as associative algebra).

\smallskip

For many of the results the fact that $k$ is a field is not essential. So we fix a commutative
ring $k$ with $1$, an associative $k$-algebra $K$, which we do not assume to be commutative, and a
$k$-algebra endomorphism $\alpha$ of $K$, and we consider the Ore extension $B=K[x,\alpha]$,
namely the algebra generated by $K$ and $x$ subject to the relations
\[
x\lambda=\alpha(\lambda)x\quad\text{for all } \lambda\in K.
\]
Let $f=x^n+\sum_{i=1}^{n}\lambda_ix^{n-i}$  be a monic polynomial of degree $n\ge 2$, where each
coefficient $\lambda_i\in K$ satisfies $\alpha(\lambda_i)=\lambda_i$ and $\lambda_i\lambda =
\alpha^i(\lambda)\lambda_i$ for every $\lambda\in K$. Sometimes we will write $f= \sum_{i=0}^n
\lambda_ix^{n-i}$, assuming that $\lambda_0=1$. For instance, the above conditions hold if the
$\lambda_i$'s are in the center of $K$, $\alpha(\lambda_i) = \lambda_i$ for all $i$ and
$\alpha^i=\id$ for all $i$ such that $\lambda_i\neq 0$. Finally let $A=B/\cl{f}$. We call $A$ the
{\em monogenic extension} of $K$ associated with $\alpha$ and $f$. Notice that under these
assumptions,
\[
fx=xf\quad\text{and}\quad f\lambda=\alpha^n(\lambda)f\quad \text{for all } \lambda\in K,
\]
and so $fB \subseteq Bf = \cl{f}$.

\smallskip

Motivated by the problem mentioned above, we study the Hochschild cohomology ring $\H_K^*(A)$ of
$A$ with coefficients in $A$, relative to $K$. If $K$ is a separable $k$-algebra then $\H_K^*(A) =
\H^*(A)$ and so, our results apply to this case. As we said above, based on this study we were
able to completely compute the cohomology ring of all the rank one Hopf algebras, which it was our
original objective. Our methods extend to the noncommutative setting those of~\cite{B}.

\smallskip

The paper is organized as follows: Section~1 is devoted to establish some notations and basic
results. In Section~2, we obtain a small resolution of $A$ as an $A$-bimodule relative to $K$ and
we build comparison maps between this resolution and the normalized canonical one. Sections~3 and
4 are the core of the paper. In the first part of Section~3 we use these results to obtain a small
cochain complex $C_S(A,M)$ given the (relative to $K$) Hochschild cohomology $\H_K^*(A,M)$, of $A$
with coefficients in $M$, for each $A$-bimodule $M$. When $M = A$ we will write $C_S(A)$ instead
of $C_S(A,A)$ and $\HH_K^*(A)$ instead of $\H_K*(A,A)$. Moreover, we obtain a map $C_S(A)\times
C_S(A)\to C_S(A)$ inducing the cup product in $\HH_K^*(A)$. Then, in Subsection~3.2, applying
these results we compute the cohomology ring of $A$, under suitable hypothesis. In Section 4 we
apply the results obtained in Section~3 to study closely the cohomology of $A =
K[x,\alpha]/\langle f \rangle$, when $K = k[G]$ is the group algebra of a finite group $G$ and
$\alpha$ is the automorphism defined by $\alpha(g) = \chi(g) g$, where $\chi\colon G\to
k^{\times}$ is a character satisfying suitable hypothesis. Finally, in Section~5 we solved the
problem possed in \cite{B-W}.

\smallskip

\noindent\bf Acknowlegetment:\rm\enspace We thank Graciela Carboni for a carefully reading of 
a preliminary version and useful remarks.

\section{Preliminaries}
In this section we fix the general terminology and notation used in the following, and establish
some basic formulas.

\medskip

Let $K$, $\alpha$, $B$, $f$ and $A$ be as in he introduction. We remark that, from the condition
$Bf = fB$ and the fact that $f$ monic, it follows that $\{1,x,\dots, x^{n-1} \}$ is a left
$K$-basis of the algebra $A$. More precisely, given $P\in B$, there exist unique $\ra{P}$ and
$\stackrel{\dots}{P}$ in $B$ such that
\[
P=\ra{P}f+\stackrel{\dots}{P}\quad\text{and}\quad \stackrel{\dots}{P}=0  \text{ or } \deg
\stackrel{\dots}{P}<n.
\]
In this paper, unadorned tensor product $\ot$ means $\ot_K$ and all the maps are $k$-linear and
all $k$-bimodules are symmetric. Given a $K$-bimodule $M$, we let $M\ot$ denote the quotient
$M/[M,K]$, where $[M,K]$ is the $k$-module generated by the commutators $m\lambda - \lambda m$
with $\lambda\in K$ and $m\in M$. Given a $k$-algebra extension $C/K$, let $C^2_{\alpha^r}
:=C_{\alpha^r}\ot C$, where $C_{\alpha^r}$ is $C$ endowed with the regular left $C$-module
structure and with the right $K$-module structure twisted by $\alpha^r$, namely, if $c\in
C_{\alpha^r}$ and $\lambda \in K$, then $c\cdot \lambda = c \alpha^r(\lambda)$. We define
\[
\frac{T}{Tx}: B\to B^2_{\alpha}
\]
as the unique $K$-derivation such that $\frac{Tx}{Tx}=1\ot 1$. Notice that
\[
\frac{Tx^i}{Tx} = \sum_{\ell=0}^{i-1}x^{\ell}\ot x^{i-\ell-1}.
\]
Composing with the canonical projection $B^2_{\alpha}\to A^2_{\alpha}$ we also obtain a
well-defined derivation $\frac{T}{Tx}\colon B\to A^2_{\alpha}$.

\begin{lem}\label{lema 1.1} On $A^2_{\alpha}$ the following equality holds for all $0\leq
i\leq n-1$:
\[
\frac{T(fx^i)}{Tx}= x^i\frac{Tf}{Tx} =\frac{Tf}{Tx}x^i.
\]
\end{lem}

\begin{proof} Since $\frac{T}{Tx}$ is a derivation,
\[
\frac{T(fx^i)}{Tx}= \frac{Tf}{Tx} x^i + f\frac{T x^i}{Tx} = \frac{Tf}{Tx} x^i,
\]
but also
\[
\frac{T(fx^i)}{Tx}= \frac{T(x^if)}{Tx} = x^i\frac{Tf}{Tx},
\]
which finish the proof.
\end{proof}

\section{The resolution}
Let $K$, $\alpha$ and $f$ be as in the introduction and let $A$ be the monogenic extension of $K$
associated with $\alpha$ and $f$. Let $\Upsilon$ be the family of all $A$-bimodule epimorphisms
which split as $K$-bimodule maps. It is easy to see that the $A_{\alpha^r}$'s are
$\Upsilon$-projective. The aim of this section is to obtain a $\Upsilon$-projective resolution
$C_S'(A)$, of $A$ as an $A$-bimodule, smaller than the normalized canonical one. We also build
comparison maps between $C_S'(A)$ and the normalized canonical resolution. We consider the
following complex
\[
\spreaddiagramcolumns{-0.3pc} \wt{C}'_S(A)=  \xymatrix@C-8pt{\cdots \ar[r]& A^2_{\alpha^{2n+1}}
\ar[r]^-{d'_5}& A^2_{\alpha^{2n}} \ar[r]^-{d'_4}& A^2_{\alpha^{n+1} } \ar[r]^-{d'_3}&
A^2_{\alpha^n} \ar[r]^-{d'_2}& A^2_{\alpha } \ar[r]^-{d'_1}& A^2 \ar[r]^m& A },
\]
where $m$ is the multiplication map and the differentials $d'_*$ are the \hbox{$A$-bimodule} maps
\[
d'_{2m+1}\colon  A^2_{\alpha^{mn+1}}\to A^2_{\alpha^{mn}} \quad \text{and} \quad d'_{2m}\colon
A^2_{\alpha^{mn}}\to A^2_{\alpha^{(m-1)n+1}},
\]
defined by
\begin{align*}
& d'_{2m+1}(1\ot 1) = x\ot 1-1\ot x,\\
& d'_{2m}(1\ot 1) = \frac{Tf}{Tx}= \sum_{i=1}^{n} \lambda_{n-i}\sum_{\ell=0}^{i-1}x^{\ell}\ot
x^{i-\ell-1}.
\end{align*}
We remark that, since we are tensoring over $K$ and not over $k$, the twisting by powers of
$\alpha$ in the modules is necessary for the well-definition of the maps.

\begin{teo}\label{teorema 2.1} $\wt{C}'_S(A)$ is contractible as an $(A,K)$-bimodule complex.
A contracting homotopy is given by the maps
\[
\sigma_0\colon  A\to A\ot A, \!\!\quad \!\!\sigma_{2m+1}\colon  A^2_{\alpha^{mn}}\to
A^2_{\alpha^{mn+1}}\!\! \quad \!\!\hbox{and} \!\!\quad\!\! \sigma_{2m}\colon
A^2_{\alpha^{(m-1)n+1}}\to A^2_{\alpha^{mn}},
\]
defined by
\begin{align*}
\allowdisplaybreaks
\sigma_0(a) &=  a\ot 1, \\
\sigma_{2m+1}(a\ot x^i) &= - a\frac{T x^i}{Tx},\\
\sigma_{2m}(a\ot x^i)&  = \begin{cases} 0&\text{if $i<n-1$,}\\
a\ot 1&\text{if $i=n-1$,}
\end{cases}
\end{align*}
Consequently, the complex
\[
C_S'(A) = \xymatrix@C-8pt{\cdots\ar[r]& A^2_{\alpha^{2n+1}} \ar[r]^-{d'_5}& A^2_{\alpha^{2n}}
\ar[r]^-{d'_4}& A^2_{\alpha^{n+1}}\ar[r]^-{d'_3}& A^2_{\alpha^n}\ar[r]^-{d'_2}& A^2_{\alpha}
\ar[r]^-{d'_1}& A^2 }
\]
is a $\Upsilon$-projective resolution of $A$.
\end{teo}

\begin{proof} We leave to the reader to check that these maps are well-defined. Let us check
that $\sigma$ is a contracting homotopy. To begin, it is clear that $m\sigma_0=\id_A$ and
$\sigma_0m(1\ot x^i) = x^i\ot 1$. Moreover
\begin{align*}
\allowdisplaybreaks
d'_{2m+1}\sigma_{2m+1}(1\ot x^i) &= d'_{2m+1}\left(-\frac{T x^i}{T x}\right)\\
&= - \sum\limits_{\ell= 0}^{i-1}(x^{\ell+1}\ot x^{i-\ell-1} - x^{\ell}\ot x^{i-\ell})\\
&= 1\ot x^i -x^i\ot 1.
\end{align*}
So, in particular, $(\sigma_0m+d'_1\sigma_1)(1\ot x^i) = 1\ot x^i$. Making the computation we
obtain
\begin{align*}
\sigma_{2m-1}d'_{2m-1}(1\ot x^{n-1})& =\sigma_{2m-1}(x\ot x^{n-1} - 1\ot \stackrel{\dots\dots}{x^{n}})\\
&=\sigma_{2m-1}(x\ot x^{n-1} - 1\ot (x^n - f))\\
&= -x\frac{T x^{n-1}}{Tx} +\frac{T x^n}{Tx}-\frac{T f}{Tx}\\
&=1\ot x^{n-1}-\frac{T f}{Tx}
\end{align*}
and
\begin{align*}
\allowdisplaybreaks
\sigma_{2m-1}d'_{2m-1}(1\ot x^i) &=\sigma_{2m-1}(x\ot x^i - 1\ot x^{i+1})\\
&=\frac{T x^{i+1}}{Tx} - x\frac{T x^i}{Tx}\\
&=1\ot x^i,
\end{align*}
for $i<n-1$. Besides,
\[
d'_{2m}\sigma_{2m}(1\ot x^{n-1}) = \frac{Tf}{Tx}\quad\text{and}\quad d'_{2m}\sigma_{2m}(1\ot
x^i)=0,\, \text{ for $i<n-1$}.
\]
So, $d'_{2m}\sigma_{2m}+ \sigma_{2m-1}d'_{2m-1}=\id$. To finish the proof, it remains to check
that
\[
d'_{2m+1}\sigma_{2m+1} + \sigma_{2m}d'_{2m} = \id.
\]
But, by Lemma~\ref{lema 1.1} and the fact that $\sigma_{2m}$ is left $A$-linear,
\[
\sigma_{2m}d'_{2m}(1\ot x^i)= \sigma_{2m} \left(x^i\frac{Tf}{Tx} \right) = x^i\ot 1 = 1\ot x^i -
d'_{2m+1}\sigma_{2m+1}(1\ot x^i),
\]
as desired.
\end{proof}

\subsection{Comparison maps} From now on we will use the standard notations
$$
\ra{A} = A/K,\quad \ra{A}^{\ot^0} = K\quad\text{and}\quad \ra{A}^{\ot^n} = \ra{A}\ot\cdots\ot\ra{A}
\,\,\,\text{ ($n$ times)},
$$
for $n\ge 1$.

\begin{prop}\label{prop:2.2} The family of $A$-bimodule maps
\[
\psi'_{2m}\colon A\ot \ra{A}^{\ot^{2m}}\ot A\to A^2_{\alpha^{mn}} \quad \text{and}\quad
\psi'_{2m+1}\colon A\ot \ra{A}^{\ot^{2m+1}}\ot A\to A^2_{\alpha^{mn+1}},
\]
recursively defined by $\psi'_0=\id$ and
\[
\psi'_{n+1}(1\ot a_1\ot \cdots\ot a_{n+1}\ot 1) =\sigma_{n+1}\psi'_n b'_{n+1}(1\ot a_1\ot\cdots
\ot a_{n+1}\ot 1),
\]
is an homotopy equivalence, from $(A\ot \ra{A}^{\ot^*}\ot A,b')$ to $C_S'(A)$, with homotopy
inverse given also recursively by $\phi'_0=\id$ and
$$
\phi'_{n+1}(1\ot 1)=\zeta_{n+1} \phi'_n d'_{n+1}(1\ot 1),
$$
where $\zeta_{n+1}(a_0\ot\cdots \ot a_{n+1})= (-1)^{n+1} \ot a_0\ot\cdots \ot a_{n+1}\ot 1$.
\end{prop}

\begin{proof} It follows by a standard argument in homological algebra.
\end{proof}

Recall that $\stackrel{\dots}{P}$ and $\ra{P}$ are defined as the unique polynomials such that
\[
\]
and $\stackrel{\dots}{P}= 0$ or $\deg P< n$.

\begin{teo}\label{teorema 2.3} The explicit formulas for $\phi'$ and $\psi'$ are
\begin{align*}
& \phi'_0(1\ot 1)= 1\ot 1,\\
& \phi'_1(1\ot 1)= 1\ot x\ot 1,\\
&\phi'_{2m}(1\ot 1) = \sum_{\mathbf{i}\in \I_m} \boldsymbol{\lambda}_{\mathbf{n}- \mathbf{i}}
\sum_{\boldsymbol{\ell}\in \J_{\mathbf{i}}}
x^{|\mathbf{i}-\boldsymbol{\ell}|-m} \ot \wt{\x}^{\boldsymbol{\ell}_{m,1}}\ot 1,\\
&\phi'_{2m+1}(1\ot 1) = \sum_{\mathbf{i}\in \I_m} \boldsymbol{\lambda}_{\mathbf{n}- \mathbf{i}}
\sum_{\boldsymbol{\ell}\in \J_{\mathbf{i}}}
x^{|\mathbf{i}-\boldsymbol{\ell}|-m} \ot \wt{\x}^{\boldsymbol{\ell}_{m,1}}\ot x\ot 1,\\
&\psi'_{2m}(1\ot\x^{\i_{1,2m}}\ot 1)= \ra{x^{i_1+i_2}}\, \ra{x^{i_3+i_4}}\cdots
\ra{x^{i_{2m-1}+i_{2m}}}\ot 1,\\
&\psi'_{2m+1}(1\ot\x^{\i_{1,2m+1}}\ot 1)= \ra{x^{i_1+i_2}}\, \ra{x^{i_3+i_4}} \cdots
\ra{x^{i_{2m-1}+i_{2m}}}\frac{T (x^{i_{2m+1}})}{Tx} ,
\end{align*}
where

\begin{itemize}

\smallskip

\item  $\I_m=\{(i_1,\dots,i_m)\in \mathbb{Z}^m: 1\le i_j\le n  \text{ for all $j$}\}$,

\smallskip

\item $\J_{\mathbf{i}}=\{(l_1,\dots,l_m)\in \mathbb{Z}^m: 1\le l_j<i_j \text{ for all
$j$}\}$,

\smallskip

\item $\boldsymbol{\lambda}_{\mathbf{n}- \mathbf{i}} =\lambda_{n-i_1}\cdots \lambda_{n-i_m}$,

\smallskip

\item $\wt{\x}^{\boldsymbol{\ell}_{m,1}} = x\ot x^{\ell_m} \ot \cdots\ot x \ot x^{\ell_1}$,

\smallskip

\item $|\mathbf{i}-\boldsymbol{\ell}|=\sum_{j=1}^m (i_j -\ell_j)$.

\smallskip

\item $\x^{\i_{1r}} = x^{i_1}\ot \cdots\ot x^{i_r}$,

\end{itemize}
\end{teo}

\begin{proof} Both, the formulas for $\phi'$ and $\psi'$ can be checked by induction on the
degree. The computations for $\phi'$ are easy and straightforward.  It is clear that the equality
for $\psi'$ is true for $2m=0$, since $\psi'_0(1\ot 1)=1\ot 1$. To abbreviate, given $1\le j\le l$
we write $\x^{\i_{jl}} = x^{i_j}\ot \cdots\ot x^{i_l}$. Assume the formula for $2m$. The recursive
definition of $\psi'_{2m+1}$ gives
\begin{align*}
\psi'(1\ot\x^{\i_{1,2m+1}}\ot 1)&=\sigma\psi'b'(1\ot\x^{\i_{1,2m+1}}\ot 1)\\
&=\sigma\psi'\bigl(b'(1\ot \x^{\i_{1,2m+1}})\ot 1 - 1\ot\x^{\i_{1,2m+1}}\bigr)\\
&=\sigma(\cdots\ot 1) -\sigma(\ra{x^{i_1+i_2}} \cdots \ra{x^{i_{2m-1}+i_{2m}}}\ot x^{i_{2m+1}}) \\
&=\ra{x^{i_1+i_1}} \cdots \ra{x^{i_{2m-1}+i_{2m}}} \frac{T(x^{i_{2m+1}})}{Tx}.
\end{align*}
Assume now the formula for $2m-1$. By the recursive definition of $\psi'_{2m+1}$, we have
\begin{align*}
\psi'(1\ot\x^{\i_{1,2m}}\ot 1) & =\sigma\psi'b'(1\ot\x^{\i_{1,2m}}\ot 1)\\
&=\sigma\psi'\bigl(b'(1\ot \x^{\i_{1,2m}})\ot 1 + 1\ot\x^{\i_{1,2m}}\bigr)\\
&=\sigma\left(\ra{x^{i_1+i_2}} \cdots \ra{x^{i_{2m-3}+i_{2m-2}}}
\frac{T(x^{i_{2m-1}})}{Tx}x^{i_{2m}}\right).
\end{align*}
We note now that in $B = k[x,\alpha]$,
\begin{align*}
\frac{T(x^{i_{2m-1}})}{Tx}x^{i_{2m}} & = \frac{T(x^{i_{2m-1}+i_{2m}})}{Tx} -
x^{i_{2m-1}}\frac{T(x^{i_{2m}})}{Tx}\\
& = \frac{T(\ra{x^{i_{2m-1}+i_{2m}}}f)}{Tx}+
\frac{T(\stackrel{\dots\dots\dots\dots\dots\dots\dots}{x^{i_{2m-1}+i_{2m}}})}{Tx}
 - x^{i_{2m-1}}\frac{T(x^{i_{2m}})}{Tx},
\end{align*}
and this equality in $A$ becomes
$$
\frac{T(x^{i_{2m-1}})}{Tx}x^{i_{2m}} = \ra{x^{i_{2m-1}+i_{2m}}}\frac{T(f)}{Tx}+
\frac{T(\stackrel{\dots\dots\dots\dots\dots\dots\dots}{x^{i_{2m-1}+i_{2m}}})}{Tx} -
x^{i_{2m-1}}\frac{T(x^{i_{2m}})}{Tx},
$$
since $\frac{T}{Tx}$ is a derivation and $f$ vanishes in $A$. So, $\sigma_{2m}$,
\begin{align*}
\psi'(1\ot \x^{\i_{1,2m}}\ot 1) & = \sigma\left(\ra{x^{i_1+i_2}}\cdots\ra{x^{i_{2m-1}+i_{2m}}}
\frac{T(f)}{Tx}\right)\\
& = \ra{x^{i_1+i_2}} \cdots \ra{x^{i_{2m-1}+i_{2m}}}\ot 1,
\end{align*}
as desired.
\end{proof}

\section{Hochschild cohomology of monogenic extensions}
Let $K$, $\alpha$, $f$ and $A$ be as in he introduction. As usual, we let $A^e$ denote the
enveloping algebra $A\ot_k A^{\op}$. Given an $A$-bimodule $M$, we let $M^{\alpha^r}$ denote the
$k$-submodule
\[
M^{\alpha^r}= \{\m\in M:\m\lambda=\alpha^r(\lambda)\m\text{ for all }\lambda\in K\}\subseteq M
\]
\begin{teo}\label{teorema 3.1} Let $M$ be an $A$-bimodule. The following facts hold:
\begin{enumerate}

\smallskip

\item The cochain complex
\[
\spreaddiagramcolumns{-1pc} \qquad C_S(A,M)= \xymatrix{ \cdots &M^{\alpha^{2n}}  \ar[l]_-{d^5}
&M^{\alpha^{n+1}} \ar[l]_-{d^4} &M^{\alpha^n}     \ar[l]_-{d^3} &M^{\alpha}       \ar[l]_-{d^2}
&M^{\alpha^0}     \ar[l]_-{d^1}, }
\]
where the coboundaries are the maps defined by
\[\qquad\quad
d^{2m+1}(\m)=x\m-\m x\!\quad\text{and}\quad\! d^{2m}(\m) =\sum_{i=1}^{n}\sum_{\ell=0}^{i-1}
\lambda_{n-i} x^{\ell}\m x^{i-\ell-1},
\]
computes $\H_K^*(A,M)$.

\smallskip

\item The maps
\begin{align*}
&\phi^*\colon (\Hom_{K^e}(\ra{A}^{\ot^*},M),b^*)\to C_S(A,M),\\
&\psi^*\colon C_S(A,M)\to (\Hom_{K^e}(\ra{A}^{\ot^*},M),b^*),
\end{align*}
defined by
\begin{align*}
\qquad\quad & \phi^0(g)= g(1),\\
& \phi^1(g)= g(x),\\
&\phi^{2m}(g) =\sum_{\mathbf{i}\in\I_m} \boldsymbol{\lambda}_{\mathbf{n}- \mathbf{i}}
\sum_{\boldsymbol{\ell}\in \J_{\mathbf i}} x^{|\mathbf{i}-\boldsymbol{\ell}|-m}
g(\wt{\x}^{\boldsymbol{\ell}_{m,1}}),\\
&\phi^{2m+1}(g)=\sum_{\mathbf{i}\in\I_m} \boldsymbol{\lambda}_{\mathbf{n}- \mathbf{i}}
\sum_{\boldsymbol{\ell}\in \J_{\mathbf i}} x^{|\mathbf{i}-\boldsymbol{\ell}|-m}
g(\wt{\x}^{\boldsymbol{\ell}_{m,1}}\ot x),\\
&\psi^{2m}(\m)(\x^{\i_{1,2m}})= \ra{x^{i_1+i_2}}\cdots \ra{x^{i_{2m-1}+i_{2m}}}\m,\\
&\psi^{2m+1}(\m)(\x^{\i_{1,2m+1}}) = \!\!\sum_{\ell=0}^{i_{2m+1}-1} \ra{x^{i_1+i_2}}\cdots
\ra{x^{i_{2m-1}+i_{2m}}}x^{\ell}\m x^{i_{2m+1}-\ell-1},
\end{align*}
where we are using the same notations as in Theorem~\ref{teorema 2.3}, are chain morphisms which
are inverse one of each other up to homotopy.
\end{enumerate}
\end{teo}

\begin{proof} For the first item, apply the functor $\Hom_{A^e}(-,M)$ to the resolution
$C_S'(A)$, and use the identification
\[
\xymatrix@R-4ex{\Hom_{A^e}(A^2_{\alpha^r},M)\ar[r]^-{\cong}&M^{\alpha^r}
\\
g\ar@{|->}[r] &g(1\ot 1).}
\]
Let $\psi_*$ and $\phi_*$ be the morphism induced by comparison maps $\psi'_*$ and $\phi'_*$,
introduced in Proposition~\ref{prop:2.2}. The second item is a straightforward consequence of
Theorem~\ref{teorema 2.3}.
\end{proof}

Recall from the introduction that when $M = A$ we write $C_S(A)$ instead of $C_S(A,A)$.

\subsection{Cup product}
In this subsection we compute the cup product of $\HH_K^*(A)$ in terms of the small complex
$C_S(A)$. Given $\m_1\in C_S^p(A)$ and $\m_2 \in C_S^q(A)$, we write $\m_1\bullet \m_2=
\phi^{p+q}(\psi^p(\m_1)\smile\psi^q(\m_2))$, where $\psi^p(\m_1)\smile\psi^q(\m_2)$ denotes the
cup product (at the level of cochains) in the canonical normalized Hochschild cochain complex.

\begin{teo}\label{teorema 3.2} Let $a_1\in C_S^p(A)$ and $a_2\in C_S^q(A)$. The following
facts hold:
\begin{enumerate}

\smallskip

\item If $p$ is even or $q$ is even, then $a_1\bullet a_2 =a_1a_2$.

\smallskip

\item If $p$ and $q$ are odd, then
$$
a_1\bullet a_2=\sum\limits_{i=2}^n\sum \limits_{j_1,j_2,j_3 \ge 0\atop j_1+j_2+j_3=i-2}
\lambda_{n-i} x^{j_1}a_1 x^{j_2}a_2x^{j_3}.
$$
\end{enumerate}
\end{teo}

\begin{proof} We use the same notations as in Theorem~\ref{teorema 3.1}. We recall that the cup
product in terms of the nomalized resolution relative to $K$ is given by
\[
(g\smile h)(a_1\ot\cdots \ot a_{p+q})= g(\a_{1p}) h(\a_{p+1,p+q})
\]
for
$$
g\in\Hom_{K^e}(\ra{A}^{\ot^p}, A)\quad\text{and}\quad h\in\Hom_{K^e}(\ra{A}^{\ot^q},A),
$$
where $\a_{1p} = a_1\ot\cdots\ot a_p$ and $\a_{p+1,p+q} = a_{p+1}\ot\cdots\ot a_{p+q}$. Assume
first that $p=2u$ and $q=2v$. We have:
\begin{align*}
a_1\bullet a_2 &= \sum_{\mathbf{i}\in\I_m}\boldsymbol{\lambda}_{\mathbf{n} \boldsymbol{-}
\mathbf{i}} \sum_{\boldsymbol{\ell}\in\J_{\mathbf{i}} }x^{|\mathbf{i}-\boldsymbol{\ell}|-m}
(\psi^p(a_1)\smile \psi^q(a_2))(\wt{\x}^{\boldsymbol{\ell}_{m,1}}) \\
&= \sum_{\mathbf{i}\in\I_m}\boldsymbol{\lambda}_{\mathbf{n} \boldsymbol{-} \mathbf{i}}
\sum_{\boldsymbol{\ell}\in\J_{\mathbf{i}} } x^{|\mathbf{i}-\boldsymbol{\ell}|-m}
a_1\ra{x^{1+\ell_m}}\cdots \ra{x^{1+\ell_{v+1}}} a_2\ra{x^{1+\ell_v}}\cdots \ra{x^{1+\ell_1}},
\end{align*}
where $m=u+v$. But $\ra{x^{\ell+1}}=0$ unless $\ell=n-1$, and in  this case $\ra{x^n}=1$. If one
wants to consider non-vanishing terms, then all $\ell_j$ must be equal to $n-1$ and this also
forces that all  $i_j$ are equal to $n$. So, the sum reduce to the single term $\lambda_0^m
a_1a_2=a_1a_2$. If $p=2u+1$ and $q=2v$, or $p=2u$ and $q=2v+1$, then a similar argument as above
shows that $a_1\bullet a_2 = a_1a_2$. Finally, if $p=2u+1$ and $q=2v+1$, then \allowdisplaybreaks
\begin{align*}
a_1\bullet a_2 &= \sum_{\mathbf{i}\in\I_m} \boldsymbol{\lambda}_{\mathbf{n}\boldsymbol{-}
\mathbf{i}}\sum_{\boldsymbol{\ell}\in\J_{\mathbf{i}}}x^{|\mathbf{i}-\boldsymbol{\ell} |-m}
(\psi^p(a_1)\smile\psi^q (a_2))(\wt{\x}^{\boldsymbol{\ell}_{m,1}})\\
&=\sum_{\mathbf{i}\in\I_m} \boldsymbol{\lambda}_{\mathbf{n} \boldsymbol{-} \mathbf{i}}
\sum_{\boldsymbol{\ell}\in\J_{\mathbf{i}}}\sum_{j=0}^{\ell_1-1}x^{|\mathbf{i}- \boldsymbol{\ell}|-m}
\Gamma_{\boldsymbol{\ell}_{m,v+2}}a_1 \Gamma_{\boldsymbol{\ell}_{v+1,2}} x^ja_2 x^{\ell_1-j-1} \\
&= \sum_{i=2}^n\lambda_{n-i}\sum_{\ell=1}^{i-1}\sum_{j=0}^{\ell-1}x^{i-\ell-1}a_1 x^ja_2 x^{\ell-j-1}\\
&=\sum\limits_{i=2}^n \sum\limits_{j_1,j_2,j_3\ge 0 \atop j_1+j_2+j_3=i-2}\lambda_{n-i}x^{j_1} a_1
x^{j_2}a_2x^{j_3},
\end{align*}
where $m=u+v+1$,
$$
\Gamma_{\boldsymbol{\ell}_{m,v+2}} = \ra{x^{1+\ell_m}}\cdots
\ra{x^{1+\ell_{v+2}}}\quad\text{and}\quad \Gamma_{\boldsymbol{\ell}_{v+1,2}} =
\ra{x^{\ell_{v+1}+1}}\cdots\ra{x^{\ell_2+1}}.
$$
This finish the proof.
\end{proof}

\subsection{Explicit computations}
Let $K$, $\alpha$ and $f$ be as in the introduction and let $A$ be the corresponding monogenic
extension. In this subsection we are going to compute the cohomology of $A$ with coefficients in
$A$, under suitable hypothesis.

\begin{teo}\label{teorema 3.3} Let $C_S^r(A)$ denote the $r$-th module of $C_S(A)$. If there
exists $\breve{\lambda}\in \Z(K)$ such that
\begin{itemize}

\smallskip

\item $\alpha^n(\breve{\lambda})=\breve{\lambda}$,

\smallskip

\item $\breve{\lambda}-\alpha^i(\breve{\lambda})$ is not a zero divisor of $K$ for $1\le i<n$,

\end{itemize}
then
\[
C_S^r(A) = \begin{cases} K^{\alpha^{mn}}&\text{if $r=2m$,}\\ K^{\alpha^{mn}}x& \text{if
$r=2m+1$.}\end{cases}
\]
\end{teo}

\begin{proof} By Theorem~\ref{teorema 3.1} we know that
$$
C_S^r(A) = \begin{cases} A^{\alpha^{mn}}&\text{if $r=2m$,}\\ A^{\alpha^{mn+1}}& \text{if
$r=2m+1$.}\end{cases}
$$
Moreover it is immediate that $a = \sum_{i=0}^{n-1} \lambda'_ix^i\in A$ satisfies $a\lambda =
\alpha^r(\lambda)a$ for all $\lambda\in K$ if and only if each $\lambda'_ix^i$ satisfies the same
condition. Hence, in order to prove the theorem it will be sufficient to check that
$\lambda'x^i\in A^{\alpha^{mn}}$ if and only if $i=0$ and $\lambda'\in K^{\alpha^{mn}}$, and that
$\lambda'x^i\in A^{\alpha^{mn+1}}$ if and only if $i=1$ and $\lambda'\in K^{\alpha^{mn}}$. If
$\lambda'x^i\in A^{\alpha^{mn}}$, then
\[
\lambda'x^i\breve{\lambda} = \alpha^{mn}(\breve{\lambda})\lambda'x^i = \breve{\lambda}
\lambda'x^i= \lambda'\breve{\lambda} x^i,
\]
since $\alpha^n(\breve{\lambda})=\breve{\lambda}$ and $\breve{\lambda}\in\Z(K)$. On the other
hand, $\lambda'x^i\breve{\lambda} = \lambda'\alpha^i(\breve{\lambda}) x^i$. So,
\[
\lambda'(\breve{\lambda}-\alpha^i(\breve{\lambda}))x^i=0,
\]
which implies $\lambda' = 0$ when $1\le i<n$, since $\breve{\lambda} - \alpha^i(\breve{ \lambda})$
is not a zero divisor of $K$ in this case. Moreover, it is clear that $\lambda'\in K\cap
A^{\alpha^{mn}}$ if and only if $\lambda'\in K^{\alpha^{mn}}$. If $\lambda'x^i\in
A^{\alpha^{mn+1}}$, then
\[
\lambda'x^i\alpha^{n-1}(\breve{\lambda}) = \alpha^{(m+1)n}(\breve{\lambda})\lambda' x^i =
\breve{\lambda} \lambda' x^i= \lambda'\breve{\lambda} x^i,
\]
and $\lambda'x^i \alpha^{n-1}(\breve{\lambda}) = \lambda'\alpha^{n-1+i}(\breve{\lambda}) x^i$. So
\[
\lambda'(\breve{\lambda}- \alpha^{i+n-1}(\breve{\lambda})) x^i = 0,
\]
which implies that $\lambda' = 0$ or $i=1$. Moreover, it is easy to check that $\lambda' x\in
A^{\alpha^{mn+1}}$ if and only if $\lambda'\in K^{\alpha^{mn}}$.
\end{proof}

\begin{teo}\label{teorema 3.4} Under the hypothesis of Theorem~\ref{teorema 3.3}, the
coboundaries maps of $C^S(A)$ are given by
\[
d^{2m+1}(\lambda )= (\alpha(\lambda)-\lambda)x \quad\text{and}\quad d^{2m+2}(\lambda
x)=-\sum_{\ell=0}^{n-1}\alpha^{\ell}(\lambda) \lambda_n.
\]
Consequently, if $\lambda_n=0$, then the even coboundary maps are zero.
\end{teo}

\begin{proof} The first assertion is immediate. Let us check the second one. Let
$\lambda x \in C_S^{2m+1}(A) = K^{\alpha^{mn}}x$. By definition
\begin{align*}
d^{2m+2}(\lambda x) & = \sum_{i=1}^n\sum_{\ell=0}^{i-1}\lambda_{n-i} x^{\ell}\lambda
x^{i-\ell}\\
& = \sum_{i=1}^n\sum_{\ell=0}^{i-1}\lambda_{n-i} \alpha^{\ell}(\lambda) x^i\\
& = \sum_{\ell=0}^{n-1} \alpha^{\ell}(\lambda) x^n + \sum_{i=1}^{n-1}\sum_{\ell=0}^{i-1}
\lambda_{n-i} \alpha^{\ell}(\lambda) x^i\\
& =-\sum_{\ell=0}^{n-1}\sum_{i=0}^{n-1}\alpha^{\ell}(\lambda)\lambda_{n-i} x^i
+ \sum_{i=1}^{n-1}\sum_{\ell=0}^{i-1}\lambda_{n-i} \alpha^{\ell}(\lambda) x^i\\
&= -\sum_{\ell=0}^{n-1} \alpha^{\ell}(\lambda)\lambda_n,
\end{align*}
where the last equality follows from Theorem~\ref{teorema 3.3}.
\end{proof}

Theorem~\ref{teorema 3.4} implies that $\alpha^n(\lambda)\lambda_n = \lambda\lambda_n$ for all
$\lambda\in K^{\alpha^{mn}}$. Indeed, this can be proved directly from the hypothesis at the
beginning of this paper and then it is true with full generality. In fact,
$$
\lambda\lambda_n = \alpha^{mn}(\lambda_n)\lambda = \lambda_n\lambda = \alpha^n(\lambda) \lambda_n,
$$
where the first equality follows from the fact that $\lambda\in K^{\alpha^{mn}}$.

\begin{coro}\label{corolario 3.5} Under the hypothesis of Theorem~\ref{teorema 3.3},
\begin{align*}
& \HH_K^0(A)=\ker(\alpha-\id)\cap \Z(K),\\
& \HH_K^{2m+1}(A)= \frac{\left\{\lambda x\in K^{\alpha^{mn}}x:\sum_{\ell=0}^{n-1}\alpha^{\ell}
(\lambda)\lambda_n=0\right\}}{(\alpha-\id)(K^{\alpha^{mn}})x},\\
& \HH_K^{2m+2}(A)= \frac{\ker(\alpha -\id)\cap K^{\alpha^{(m+1)n}}}{\left\{\sum_{\ell=0}^{n-1}
\alpha^{\ell}(\lambda)\lambda_n:\lambda\in K^{\alpha^{mn}}\right\}}.
\end{align*}
\end{coro}

\subsubsection{Cup product} It is easy to refine Corollary~\ref{corolario 3.5} by describing
the cup product in the Hochschild cohomology of the extension $A/K$. By item~(1) of
Theorem~\ref{teorema 3.2} we know that the product of two homogeneous elements is induced by the
multiplication map in $K$, whenever at least one of them have even degree. On the other hand, if
$\lambda x\in \HH_K^{2m+1}(A)$ and $\lambda' x\in \HH_K^{2m'+1}(A)$, then, by item~(2) of
Theorem~\ref{teorema 3.2}, \allowdisplaybreaks
\begin{align*}
\lambda x\bullet \lambda' x &= \sum\limits_{i=2}^n \sum\limits_{j_1,j_2,j_3\ge 0 \atop
j_1+j_2+j_3=i-2} \lambda_{n-i} x^{j_1}\lambda x x^{j_2}\lambda' x x^{j_3}\\
& = \sum\limits_{i=2}^n \sum\limits_{j_1,j_2,j_3\ge 0 \atop j_1+j_2+j_3=i-2} \lambda_{n-i}
\alpha^{j_1}(\lambda)\alpha^{j_1+j_2+1}(\lambda')x^{i}\\
& = \sum\limits_{i=2}^{n-1} \!\! \sum\limits_{j_1,j_2,j_3\ge 0 \atop j_1+j_2+j_3=i-2} \!\!
\lambda_{n-i} \alpha^{j_1}(\lambda)\alpha^{j_1+j_2+1}(\lambda')x^i\\
& - \sum\limits_{i=0}^{n-1} \!\! \sum\limits_{j_1,j_2,j_3\ge 0 \atop j_1+j_2+j_3=n-2} \!\!
\alpha^{j_1}(\lambda)\alpha^{j_1+j_2+1}(\lambda')\lambda_{n-i}x^i\\
&=-\sum\limits_{j_1,j_2,j_3\ge 0\atop j_1+j_2+j_3=n-2}\alpha^{j_1}(\lambda)\alpha^{j_1+j_2+1}
(\lambda')\lambda_n,
\end{align*}
where the last equality follows from Theorem~\ref{teorema 3.3}.

\subsubsection{Gerstenhaber structure} The goal of this paragraph is to compute the full structure 
of $\HH_K^*(A)$ as Gerstenhaber algebra, namely, to compute the Gerstenhaber bracket on $\HH_K^*(A)$, 
when $A$ is an algebra satisfying hypothesis of Theorem~\ref{teorema 3.3}. We recall first the 
definition of the Lie bracket on $\HH_K^*(A)$ introduced in \cite{G}.

\begin{defi}\label{definicion 3.6}
[Gerstenhaber] Let $C$ be a ring and $V$ a $C$-bimodule. Let $f\in\Hom_{C^e}(V^{\ot^r},V)$ and
$g\in\Hom_{C^e}(V^{\ot^{r'}},V)$. The composition into the $j$-th place of $f$ and $g$, is the map
$$
f\circ_j g\in \Hom_{C^e}(V^{\ot^{r+r'-1}},V),
$$
defined by
\[
f\circ_j g (\v_{1,r+r'-1})=f(\v_{1,j-1}\ot g(\v_{j,j+r'-1})\ot \v_{j+r',r+r'-1}),
\]
where $\v_{h,l} = v_h\ot\cdots\ot v_l$ for $h\le l$. The composition product $f\circ g$ and the
bracket $[f,g]$ are defined by
\begin{align*}
& f\circ g = \sum_{j=1}^r(-1)^{(j+1)(r'+1)}f\circ_j g\\
& [f,g] = f\circ g-(-1)^{(r+1)(r'+1)}g\circ f.
\end{align*}
\end{defi}

Since our comparison maps are between $C^S(A)$ and the relative to $K$ normalized complex
$\bigl(\Hom_{K^e}(\ra{A}^{\ot^*},A),b^*\bigr)$, to compute the Gerstenhaber bracket we need to
identify $\Hom_{K^e}(\ra{A}^{\ot^*},A)$ with the $K^e$-submodule of $\Hom_{K^e}(A^{\ot^*},A)$
consisting of all the functions $f$ vanishing on all the simple tensors $a_1\ot\cdots\ot
1\ot\cdots\ot a_p$. It is well-known that in this way one obtain a subcomplex of
$\bigl(\Hom_{K^e}(A^{\ot^*},A),b^*\bigr)$ and that the canonical inclusion is a quasi-isomorphism.
Also, it is clear that this subcomplex is, in fact, a Lie subalgebra. In order to establish the
main result of this paragraph we need first to introduce some notations

\begin{defi}\label{definicion 3.7} For $a\in C_S^{r}(A)$, $a'\in C_S^{r'}(A)$ and  $j=1,\dots,r$,
we define
\begin{align*}
& a\circ^s_j a' = \phi^{r+r'-1}\left(\psi^{r}(a)\circ_j \psi^{r'}(a') \right),\\
& a\circ^s a'= \sum_{j=1}^r(-1)^{(j+1)(r'+1)} a\circ^s_j a',\\
& [a,a']_s= a\circ^s a' - (-1)^{(r+1)(r'+1)} a'\circ^s a.
\end{align*}
\end{defi}

By construction $[-,-]_s$ induce the Gerstenhaber bracket in $\HH_K^*(A)$.

\begin{teo}\label{teorema 3.8} Let $\lambda\in K^{\alpha^{mn}}$ and $\mu\in K^{\alpha^{m'n}}$. Assume that we are
in the hypothesis of Theorem~\ref{teorema 3.3}. We have:
\begin{align*}
& [\lambda,\mu]_s = 0,\\
& [\lambda,\mu x]_s = \sum_{h=0}^{mn-1} \alpha^h(\mu)\lambda,\\
& [\lambda x,\mu x]_s = \sum_{h=0}^{mn} \alpha^h(\mu)\lambda - \sum_{h=0}^{m'n}
\alpha^h(\lambda)\mu.
\end{align*}

\end{teo}

To prove this theorem we are going to use the following lemmas:

\begin{lem}\label{lema 3.9} Let $\mu\in K^{\alpha^{m'n}}$. The following equalities are true:
\begin{align*}
& \psi^{2m'}(\mu)(\wt{\x}^{\boldsymbol{\ell}_{m',1}}) = \begin{cases} \mu & \hbox{if $\ell_1 =
\dots = \ell_{m'} = n-1$,}\\ 0  & \hbox{otherwise,}\end{cases}\\
& \psi^{2m'}(\mu)(x^{\ell_{m'}}\ot \wt{\x}^{\boldsymbol{\ell}_{m'-1,1}}\ot x) = \begin{cases} \mu
&
\hbox{if $\ell_1 = \dots = \ell_{m'} = n-1$,}\\ 0  & \hbox{otherwise,}\end{cases}\\
& \psi^{2m'+1}(\mu x)(\wt{\x}^{\boldsymbol{\ell}_{m',1}}\ot x) = \begin{cases} \mu x& \hbox{if
$\ell_1 = \dots = \ell_{m'} = n-1$,}\\ 0  & \hbox{otherwise,}\end{cases}\\
& \psi^{2m'+1}(\mu x)(x^{\ell_{m'+1}}\ot\wt{\x}^{\boldsymbol{\ell}_{m',1}}) = \begin{cases}
 \delta_{l_1}(\mu) x^{l_1}& \hbox{if $\ell_2 = \dots = \ell_{m'+1} = n-1$,}\\0 & \hbox{otherwise,}\end{cases}
\end{align*}
where $\delta_{l_1}(\mu) = \sum_{h=0}^{\ell_1-1} \alpha^h(\mu)$.
\end{lem}

\begin{proof} Theses equalities follow by a direct computation, using the formulas for $\psi$
obtained in Theorem~\ref{teorema 3.4}.
\end{proof}

\begin{lem}\label{lema 3.10} Let $\lambda\in K^{\alpha^{mn}}$ and $\mu\in K^{\alpha^{m'n}}$. Let $\delta(\mu) =
\sum_{h=0}^{n-2} \alpha^h(\mu)$. The following equalities are true:

\begin{enumerate}

\item $\lambda\circ^s_j \mu = 0$,

\item $\lambda x\circ^s_j \mu = 0$,

\item If $j$ is odd, then $\lambda \circ^s_j \mu x = \alpha^{(j-1)n/2}(\mu)\lambda$,

\item If $j$ is even, then $\lambda \circ^s_j \mu x = \alpha^{1+(j-2)n/2}(\delta(\mu))\lambda$,

\item If $j$ is odd, then $\lambda x\circ^s_j \mu x = \alpha^{(j-1)n/2}(\mu)\lambda x$,

\item If $j$ is even, then $\lambda x\circ^s_j \mu x = \alpha^{1+(j-2)n/2}(\delta(\mu))\lambda x$.

\end{enumerate}

\end{lem}

\begin{proof} All the equalities follow in a similar way. We prove the last one. By definition
and the formula for $\phi^{2m+2m'+1}$ obtained in Theorem~\ref{teorema 3.4}, we have
\[
\lambda x\circ_j\mu x= \sum_{\stackrel{\mathbf{i}\in\I_u}{\boldsymbol{\ell}\in\J_{\mathbf{i}}}}
\boldsymbol{\lambda}_{\mathbf{n}\boldsymbol{-} \mathbf{i}} x^{|\mathbf{i}-\boldsymbol{\ell}|-u}
(\psi^{2m+1}(\lambda x)\circ_j \psi^{2m'+1}(\mu x))(\wt{\x}^{\boldsymbol{\ell}_{u,1}}\ot x)
\]
where $u = m+m'$. By the last equality in Lemma~\ref{lema 3.10}, we know that
\[
\lambda x\circ_j\mu x= \sum_{\mathbf{i}\in\I_m} \sum_{\boldsymbol{\ell}\in\J_{\mathbf{i}}}
\boldsymbol{\lambda}_{\mathbf{n}\boldsymbol{-} \mathbf{i}} x^{|\mathbf{i}-\boldsymbol{\ell}|-m}
\psi^{2m+1}(\lambda x)\bigl(\wt{\x}^{(\boldsymbol{\ell}_{m,1,j})}\bigr),
\]
where
$$
\wt{\x}^{(\boldsymbol{\ell}_{m,1,j})} = \wt{\x}^{\boldsymbol{\ell}_{m,m-(j-4)/2}}\ot x\ot
\delta_{\ell_{m-(j-2)/2}}(\mu) x^{\ell_{m-(j-2)/2}}\ot\wt{\x}^{\boldsymbol{\ell}_{m-j/2,1}}\ot x.
$$
But, by the formula for $\psi^{2m+1}$ obtained in Theorem~\ref{teorema 3.4},
$$
\psi^{2m+1}(\lambda x)\bigl(\wt{\x}^{(\boldsymbol{\ell}_{m,1,j})}\bigr)=\begin{cases} \alpha^{1+(j-2)n/2}(\delta(\mu))\lambda x &\text{if $l_h=n-1$ for all $h$},\\
0 &\text{otherwise}.
\end{cases}
$$
So, $\lambda x\circ^s_j \mu x = \alpha^{1+(j-2)n/2}(\delta(\mu))\lambda x$, as desired.
\end{proof}

\subsubsection*{Proof of Theorem~\ref{teorema 3.8}} The first equality follows immediately from
items~(1) and~(2) of Lemma~\ref{lema 3.10}. We check the third one and left the second one to the
reader. By definition and items~(5) and~(6) of Lemma~\ref{lema 3.10}
\begin{align*}
\lambda x \circ^s \mu x & = \sum_{h=0}^m \lambda x\circ^s_{2h+1} \mu x +
\sum_{h=1}^m \lambda x\circ^s_{2h} \mu x\\
& = \sum_{h=0}^m \alpha^{hn}(\mu) \lambda x + \sum_{h=0}^{m-1} \alpha^{hn}(\delta(\mu)) \lambda x\\
& = \sum_{h=0}^{mn} \alpha^h(\mu) \lambda x.
\end{align*}
Similarly
$$
\mu x \circ^s \lambda x =  \sum_{h=0}^{m'n} \alpha^h(\lambda)\mu x.
$$
The third equality follows from these facts.\qed

\begin{coro}\label{corolario 3.11} Under the hypothesis of Theorem \ref{teorema 3.3},
$$
[\HH_K^{\even},\HH_K^{\even}]=0.
$$
In particular $[\HH_K^2(A),\HH_K^2(A)]=0$ and so, if $K$ is separable, then  every infinitesimal deformation of $A$ may be
completed into a formal deformation.
\end{coro}

\begin{rem}\label{remark 3.12} Note that the results in Theorems~\ref{teorema 3.3} and
\ref{teorema 3.4}, Corollary~\ref{corolario 3.5} and the cup product and Gerstenhaber bracket
do not depend on $f$, with the exception of its degree $n$ and its independent term $\lambda_n$.
\end{rem}

\begin{coro}\label{corolario 3.13} Let $C_n = \langle w \rangle$ be the cyclic group with $n$
elements. Assume that the hypothesis of Theorem~\ref{teorema 3.3} are satisfied. If $\lambda_n$ is
invertible and $\alpha^n = \id$, then
$$
\HH_K^*(A) = H^*(\mathds{Z}_n,\Z(K)(\alpha)),
$$
where $\Z(K)(\alpha)$ denotes $\Z(K)$ endowed with the structure of $C_n$-module given by $w\cdot
\lambda = \alpha(\lambda)$.
\end{coro}

\begin{teo}\label{teorema 3.14} Suppose that $k$ is a field and that the hypothesis of
Theorem~\ref{teorema 3.3} are satisfied. If $\alpha$ is a diagonalizable epimorphism, then
\begin{align*}
&\HH_K^0(A) = \ker(\alpha-\id)\cap\Z(K) ,\\
&\HH_K^{2m+1}(A)=\left(\ker(\alpha-\id)\cap K^{\alpha^{mn}}\cap\Ann(n\lambda_n)\right)x,\\
&\HH_K^{2m+2}(A)=\frac{\ker(\alpha-\id)\cap K^{\alpha^{(m+1)n}}}{n\lambda_n
\left(\ker(\alpha-\id)\cap K^{\alpha^{mn}}\right)},
\end{align*}
where $\Ann(n\lambda_n) = \{\lambda\in K:n\lambda\lambda_n = 0\}= 0$.

\end{teo}

\begin{proof} Using that $\alpha$ is an epimorphism it follows easily that each $K^{\alpha^{mn}}$ is
$\alpha$-invariant. So each $K^{\alpha^{mn}}$ decompose as a direct sum
\[
K^{\alpha^{mn}}=\left(K^{\alpha^{mn}}\cap\ker(\alpha-\id)\right)\bigoplus\left(K^{\alpha^{mn}}
\cap\ker(\alpha-\id)\right)^{\perp},
\]
where $\left(K^{\alpha^{mn}}\cap\ker(\alpha-\id)\right)^{\perp}$ is the direct sum of the
eigenspaces of $\alpha$ with eigenvalue different from one. Clearly
$$
\left(K^{\alpha^{mn}}\cap\ker(\alpha-\id)\right)^{\perp}=(\alpha-\id)(K^{\alpha^{mn}}),
$$
and so
\[
\sum_{\ell=0}^{n-1}\alpha^{\ell}(\lambda)\lambda_n=0\text{ on } \left(K^{\alpha^{mn}}\cap
\ker(\alpha-\id)\right)^{\perp},
\]
since $\alpha^{n}(\lambda)\lambda_n= \lambda\lambda_n$ for all $\lambda\in  K^{\alpha^{mn}}$. The
result follows immediately from this fact and  Corollary~\ref{corolario 3.5}.
\end{proof}

\begin{rem}\label{remark 3.15} Assume that the hypothesis of Theorem~\ref{teorema 3.14} are
satisfied. From the formula for the cup product it follows that if $\lambda x\in \HH_K^{2m+1}(A)$ and
$\lambda' x\in \HH_K^{2m'+1}(A)$, then
\[
\lambda x\bullet \lambda' x = - \sum\limits_{j_1,j_2,j_3\ge 0 \atop j_1+j_2+j_3=n-2}
\alpha^{j_1}(\lambda) \alpha^{j_1+j_2+1}(\lambda')\lambda_{n} = - \binom{n}{2}\lambda\lambda'
\lambda_n,
\]
where the last equality is true since $\alpha(\lambda) = \lambda$ and $\alpha(\lambda') =
\lambda'$. Since $\HH_K^{2m+1}(A)\subseteq \Ann(n\lambda_n)$, this implies  that the product of
two elements of odd degree is zero if the characteristic of $k$ is different from~$2$, and also if
the characteristic of $k$ is $2$, but $n$ is odd or $4$ divides $n$. If the characteristic of $k$
is $2$, $n$ is even and $n/2$ is odd, then $\lambda x\bullet \lambda' x=\lambda\lambda'\lambda_n$.
\end{rem}

Next we consider another situation in which the cohomology of $A$ can be compute. The following
results are very closed to the ones valid in the commutative setting.

\begin{teo}\label{teorema 3.16} If $\alpha$ is the identity map, then
$$
C_S^r(A) = \cent(K)\oplus \cent(K)x\oplus \cdots \oplus \cent(K)x^{n-1}=\frac{\cent(K)[x]}{
\langle f\rangle},
$$
where $\cent(K)$ is the center of $K$. Moreover, the odd coboundary maps $d^{2m+1}$ of $C_S(A)$
are zero, and the even coboundary maps $d^{2m}$ are the multiplication by the derivative $f'$ of
$f$.
\end{teo}

\begin{proof} This is immediate.
\end{proof}

\begin{coro}\label{corolario 3.17} If $\alpha$ is the identity map, then
\begin{align*}
& \HH_K^0(A)= \frac{\cent(K)[x]}{\langle f\rangle},\\
& \HH_K^{2m+1}(A)=\Ann(f'),\\
& \HH_K^{2m+2}(A)=\frac{\cent(K)[x]}{\langle f,f'\rangle},
\end{align*}
where $\Ann(f')=\left\{a\in \frac{\cent(K)[x]}{\langle f\rangle}:af'=0\right\}$.
\end{coro}

\section{Cohomology of some extensions of a group algebra}
Let $k$ be a field, $K = k[G]$ the group $k$-algebra of a finite group $G$ and $\chi\colon G\to
k^{\times}$ a character. Let $\alpha\colon K\to K$ be the automorphism defined by $\alpha(g)
=\chi(g)g$ and let $f = x^n + \lambda_1 x^{n-1}+\cdots +\lambda_n\in K[x]$ be a monic polynomial
whose coefficients satisfy the hypothesis required in the introduction. Assume that there exists
$g_1\in\Z(G)$ such that $\chi(g_1)$ is a primitive $n$-th root of $1$. In this section we apply
the results obtained in Section~3 to compute the cohomological ring of $A = K[x,\alpha]/\langle
f\rangle$. Note that the hypothesis of Theorem~\ref{teorema 3.3} are fulfilled, taking $\breve
\lambda = g_1$. In particular the cohomological behavior of $A$ is independent of the polynomial
$f$, with the exception of its degree and its independent term. Since $\alpha$ is diagonalizable
Theorem~\ref{teorema 3.14} and Remark~\ref{remark 3.15} apply. In order to use the former we need
first to make some computations. By definition,
\begin{align*}
& \ker(\alpha-\id)=k[N]\text{, where }N=\ker(\chi:G\to k^{\times})\\
\intertext{and}
&K^{\alpha^{mn}}= \left\{\sum_{g\in G}\gamma_g g\in k[G]:\sum_{g\in G}\gamma_g g= \sum_{g\in G}
\gamma_g\chi^{mn}(h) h gh^{-1}\text{ $\forall\, h\in G$} \right\}.
\end{align*}
Note that $\Z(K)\cap k[N] = k[N]^G$, where $G$ acts on $k[N]$ via conjugation. By
Theorem~\ref{teorema 3.14} and the first equality, we have the following result:

\begin{teo}\label{teorema 4.1} Let $A = K[x,\alpha]/\langle f\rangle$ be as in the beginning of
this section. The Hochschild cohomology of $A/K$ is given by:
\begin{align*}
&\HH_K^0(A) = K[N]^G,\\
&\HH_K^{2m+1}(A)=\left(K[N]\cap K^{\alpha^{mn}}\cap\Ann(n\lambda_n)\right)x,\\
&\HH_K^{2m+2}(A)=\frac{K[N]\cap K^{\alpha^{(m+1)n}}}{n \left(K[N]\cap K^{\alpha^{mn}}
\right)\lambda_n}.
\end{align*}
\end{teo}

It is easy to check now that
\[
\sum_{g\in G}\gamma_g g\in K^{\alpha^r}\iff\gamma_{hgh^{-1}} = \gamma_g\chi^mn(h) \text{ for all }
h\in G.
\]
Consequently, if there exists $h\in G$ such that $\chi^r(h)\neq 1$ and $hg=gh$, then $\gamma_g=0$.
Let $\mathds{X}(r) = \{X_1,\dots,X_l\}$ be the set of the conjugation classes $X$ of $G$, which
satisfy the following property: if $h\in G$ commutes with an element of $X$, then $\chi^r(h)=1$.
It is easy to see that for each $X_j$ there exist elements $\gamma_g\in k\setminus \{0\}$, where
$g$ runs on $X_j$, satisfying
\[
\gamma_{hgh^{-1}}= \gamma_g\chi^r(h) \quad\text{for all $h\in G$.}
\]
Clearly the family $\{a_1,\dots,a_l\}$, where $a_j = \sum_{g\in X_j}\gamma_g g$, is a basis of
$K^{\alpha^r}$. Thus,
\[
K^{\alpha^r} = \bigoplus\limits_{X_j\in \mathds{X}(r)} ka_j \quad\text{and}\quad
K^{\alpha^r}\cap\ker(\alpha-\id) = \bigoplus\limits_{X_j\in \mathds{X}'(r)} ka_j,
\]
where $\mathds{X}'(r)=\{X_j\in \mathds{X}(r) : X_j\subseteq N\}$. From this discussion it follows
that $f$ satisfies the conditions required in the introduction if and only if
$$
\lambda_i\in \bigoplus\limits_{X_j\in \mathds{X}'(i)} ka_j\quad\text{for all $i$.}
$$
Moreover , by Theorem~\ref{teorema 4.1}, we have:
\begin{align*}
&\HH_K^0(A) = \bigoplus\limits_{X_j\in \mathds{X}'(0)} ka_j,\\
&\HH_K^{2m+1}(A) = \left(\Ann(n\lambda_n) \bigcap \bigoplus\limits_{X_j\in \mathds{X}'(mn)}
ka_j \right)x,\\
&\HH_K^{2m+2}(A) = \frac{\bigoplus\limits_{X_j\in \mathds{X}'((m+1)n)} ka_j}{n\left(
\bigoplus\limits_{X_j\in \mathds{X}'(mn)} ka_j\right)\lambda_n}.
\end{align*}

\begin{rem}\label{remark 4.2} Let $X$ be a conjugation class of $G$. It is easy to check that
there exists $m_0\ge 0$ such that $X\in \mathds{X}(mn)$ if and only if $m_0$ divides $m$.
Consequently, $\mathds{X}(mn)\subseteq \mathds{X}(m'n)$ whenever $m$ divides $m'$.
\end{rem}

\begin{rem}\label{remark 4.3} Let $v\ge 0$ be the order of $\chi^n$. If $m$ is congruent to
$m'$ module $v$, then $\mathds{X}'(m'n) = \mathds{X}'(mn)$ and $\mathds{X}'((m'+1)n) =
\mathds{X}'((m+1)n)$, and so
$$
\HH_K^{2m+1}(A) = \HH_K^{2m'+1}(A)\quad\text{and}\quad \HH_K^{2m+2}(A)=\HH_K^{2m'+2}(A).
$$
Hence, except for $\HH_K^0(A)$, the cohomology $\HH_K^*(A)$ is periodic of period $2v$. Moreover,
if $n\lambda_n = 0$, then
$$
\HH_K^{2v}(A) = \HH_K^0(A)\quad\text{and}\quad \HH_K^{2m+1}(A) \simeq \HH_K^{2m}(A),
$$
for all $m\ge 0$. We now are going to study the cup product in $\HH_K^*(A)$. We already know that if
$p$ or $q$ is even, then the multiplication map
$$
\HH_K^p(A)\ot_k \HH_K^q(A)\to \HH_K^{p+q}(A)
$$
is induced by the multiplication in $A$. We assert that the multiplication of two odd degree
elements is zero. In fact, since $\chi(g_1)$ has order $n$, the characteristic of $k$ is relative
prime to $n$. Consequently, by Remark~\ref{remark 3.15}, the assertion is true still when the
characteristic of $k$ is $2$.
\end{rem}


\begin{rem}\label{remark 4.4} If the characteristic of $k$ is relative prime to the order of
$G$, then $k[G]$ is a separable $k$-algebra. Hence, by \cite[Theorem 1.2]{G-S}, $\HH^*(A) =
\HH_K^*(A)$. From this and the discussion above, the computations made out in
\cite[Section~2]{B-W}, follow immediately.
\end{rem}

\begin{rem}\label{remark 4.5} Let $v$ be the order of $\chi^n$. From Remark~\ref{remark 4.3}
it follows that $\HH_K^*(A)$ is generated as a $k$-algebra by

\begin{itemize}

\item $k[N]^G$ in degree $0$,

\item an arbitrary set of generators of $\HH_K^{2m+1}(A)$ as an $\HH_K^0(A)$-module in degree
$2m+1$ for $0\le m <v$,

\item an arbitrary set of generators of $\HH_K^{2m}(A)$ as an $\HH_K^0(A)$-modu\-le in degree
$2m$ for $0< m <v$, and

\item the class of $1$ in degree $2p$.

\end{itemize}
If $n\lambda_n=0$, then the situation is simpler. In this case $\HH_K^*(A)$ is generated by

\begin{itemize}

\item $k[N]^G$ in degree $0$,

\item $x$ in degree $1$,

\item an arbitrary set of generators of $\HH_K^{2m}(A)$ as an $\HH_K^0(A)$-modu\-le in degree
$2m$ for $0<m<v$, and

\item the unit $1\in k[N]^G$ in degree $2v$.

\end{itemize}
Consequently, if $\HH_K^{2m}(A)=0$ for $0<m<v$ (for instance if $v=1$), then the algebra
$\HH_K^*(A)$ is isomorphic to $k[N]^G\ot_k k[y,x]/\langle x^2 \rangle$, where the degree of $x$ is
$1$ and the degree of $y$ is $2v$. From this it follows immediately that if the characteristic of
$k$ is relative prime to the order of $G$, then $\HH^*(A)$ is isomorphic to $k[N]^G\ot_k
k[y,x]/\langle x^2 \rangle$. When $f=x^n$ this gives Theorem~3.4 of \cite{B-W}.

\end{rem}

\subsubsection{A concrete example}
Let $G = \cl{g,h : g^u=1=h^4, hg=g^{-1}h }$. Clearly $G=\{g^jh^l : 0\leq j<u, 0\leq l<4\}$. Let $k
= \C$ and let $\chi:G\to \C^{\times}$ be the character defined by $\chi(g^jh^l)=i^l$. Consider
$f=x^2$. The hypothesis of Theorem~\ref{teorema 3.3} are satisfied with $\breve\lambda= h^2$. By
Remark~\ref{remark 4.3}, we know that the cohomology of $A$ is periodic of period $4$, and, by a
direct computation, \allowdisplaybreaks
\begin{align*}
\HH^0(A)&=\left\{\sum_{j=0}^{u-1}\gamma_j g^j : \gamma_j=\gamma_{u-j}\text{ for all }j>0
\right\},\\
\HH^1(A)& = \HH^0(A)x\\
\HH^2(A) &=\left\{\sum_{j=1}^{u-1}\gamma_j g^j : \gamma_j =-\gamma_{u-j}\text{ for all } j
\right\},\\
\HH^3(A)& = \HH^2(A)x.
\end{align*}
As a $\C$-algebra, $\HH^*(A)$ is generated by $a=g+g^{-1}$ in degree zero, $x$ in degree one,
$b=g-g^{-1}$ in degree two, and $c=1$ in degree four. The $4$-periodicity is given by the
multiplication by $c$. Notice that if $r$ is even, then $\dim_k(\HH^2(A))<\dim_k(\HH^0(A))$ and so
$\HH^*(A)$ is not isomorphic to $\HH^0(A)\otimes_k V$ for any $k$-vector space $V$.

\section{An application}
Let $k$ be a field, $G$ a finite group whose order is relative prime to the characteristic of $k$
and $\chi\colon G\to k^{\times}$ a character. Assume that there exists $g_1\in \Z(G)$ such that
$\chi(g_1)$ is a primitive $n$-th root of $1$. In particular $n$ is coprime relative to the
characteristic of $k$. In this section we compute the cohomology of $A = k[G][x,\alpha]/\langle
x^n-\xi(g_1^n-1)\rangle$ over $k$, where $\xi\in k^{\times}$ and $\alpha\in \Aut(k [G])$ is
defined by $\alpha(g) = \chi(g)g$. Recall that $\ker(\alpha - \id)= k[N]$, where $N=\ker(\chi)$.
We consider two different cases. The second one solves the problem possed in \cite{B-W}, mentioned
in the introduction.

\medskip

\noindent\bf $\boldsymbol{\chi}^{\mathbf{n}} \boldsymbol{\neq} \boldsymbol{\id}$.\rm\enspace Let
$g\in G$ such that $\chi^n(g)\neq 1$. Since
$$
g^{-1}(x^n-\xi(g_1^n-1))g = \chi^n(g)x^n - \xi(g_1^n -1),
$$
we conclude that the ideal $\cl{x^n-\xi(g_1^n-1)}$ coincides with the ideal $\cl{x^n, g_1^n-1}$.
So, the algebra $A$ is equal to $k[G/\cl{g_1^n}][x,\wt{\alpha}]/\cl{x^n}$, where $\wt{\alpha}$ is
the automorphism induced by $\alpha$. We consider now $K:=k[G/\cl{g_1^n}]$ and $f=x^n$. These data
satisfy the hypothesis of Theorem~\ref{teorema 4.1}. Moreover $K:=k[G/\cl{g_1^n}]$ is separable
over $k$ and so, $\HH^*(A) = \HH_K^*(A)$. Thus,
\begin{align*}
&\HH^0(A) = K^{G/\cl{g_1^n}},\\
&\HH^{2m+1}(A) = \left(k[N/\cl{g_1^n}]\cap K^{\wt{\alpha}^{mn}} \right)x,\\
&\HH^{2m+2}(A) = k[N/\cl{g_1^n}]\cap K^{\wt{\alpha}^{(m+1)n}}.
\end{align*}
The cup product of two homogeneous elements is zero if both of them have odd degree, and it is
induced by the multiplication in $K$, otherwise.

\medskip

\noindent\bf $\boldsymbol{\chi}^{\mathbf{n}} \boldsymbol{=} \boldsymbol{\id}$.\rm\enspace In this
case $f=x^n-\xi(g_1^n-1)$ satisfies the hypothesis required in the introduction (that is
$\alpha(\xi(g_1^n-1)) = \xi(g_1^n-1)$ and $\xi(g_1^n-1)\lambda = \alpha^n(\lambda)\xi(g_1^n-1)$).
Moreover $K = k[G]$ is separable over $k$ and so, $\HH^*(A) = \HH_K^*(A)$. Since $\chi^n  = \id$,
we have
$$
K^{\alpha^{mn}} = \Z(K)\quad\text{and so}\quad K[N]\cap K^{\alpha^{mn}} = K[N]^G.
$$
Hence, by Theorem~\ref{teorema 4.1},
\allowdisplaybreaks
\begin{align*}
&\HH^0(A)=k[N]^G,\\
&\HH^{2m+1}(A)=\left(k[N]^G\cap\Ann(g_1^n-1)\right)x,\\
&\HH^{2m+2}(A)=\frac{k[N]^G}{(g_1^n-1) k[N]^G},
\end{align*}
since $n\xi\in k^{\times}$. By Remark~\ref{remark 4.3} the cup product of two homogeneous elements
is zero if both of them have odd degree, and it is induced by the multiplication in $K$,
otherwise. Moreover, using that $\chi^n = \id$ and Theorem~\ref{teorema 3.8} it is easy to see that the Gerstenhaber bracket is induce by the map $[-,-]_s\colon C_S^*(A)\times C_S^*(A)\to C_S^*(A)$ given by 
$$
[\lambda,\mu]_s = 0,\quad [\lambda,\mu x]_s = 0\quad\text{and}\quad [\lambda x,\mu x]_s = \lambda\mu-\mu\lambda, 
$$
for $\lambda\in K^{\alpha^{mn}}$ and $\mu\in K^{\alpha^{m'n}}$.

Finally let $\wt{A} = k[G/\cl{g_1^n}][x,\wt{\alpha}]/\cl{x^n}$, where $\wt{\alpha}$ is
the automorphism induced by $\alpha$. Using the formulas obtained in Section~4, it is easy to see
that $\HH^m(A) = \HH^m(\wt{A})$, for all $m>0$.

\section{A final example}
Here we consider an example in order to show that the cochain complex introduced in
Theorem~\ref{teorema 3.1} can be used to perform explicit computations still when the hypothesis
introduced in Subsection~3.2 are not satisfied.

\smallskip

Let $\mathbb{H}$ be the skew field of the quaternions over $\mathbb{R}$. Recall that $\mathbb{H}$
is the four dimensional real algebra with basis $\{1,i,j,k\}$, unit $1$ and multiplication defined
by $i^2 = j^2 = -1$ and $ij = -ji= k$. Let $\alpha\colon \mathbb{H}\to \mathbb{H}$ be the
$\mathbb{R}$-algebra automorphism, defined by
\begin{alignat*}{2}
&\alpha(1) = 1,&&\qquad  \alpha(i) = \cos\theta\,i + \sin\theta\,j,\\
& \alpha(k) = k,&&\qquad \alpha(j) = - \sin\theta\,i + \cos\theta\,j,
\end{alignat*}
where $0<\theta<2\pi$. So $\alpha$ acts as the rotation of angle $\theta$ with axis $k$ on the
pure imaginary quatenions. Take a monic polynomial with coefficients in $\mathbb{H}$
$$
f = x^n+\lambda_1x^{n-1}+\cdots+\lambda_n
$$
and write
$$
\lambda_u = \lambda_{u1} + \lambda_{ui}i + \lambda_{uj}j + \lambda_{uk}k\qquad\text{with
$\lambda_{u1},\lambda_{ui},\lambda_{uj},\lambda_{uk}\in \mathbb{R}$.}
$$
Next we ask for the conditions in order that $f$ satisfies the hypothesis required in the
introduction. That is
$$
\alpha(\lambda_u)=\lambda_u\quad\text{ and} \quad\lambda_u\lambda = \alpha^u(\lambda)
\lambda_u\text{ for all $\lambda\in \mathbb{H}$.}
$$
By definition
$$
\alpha(\lambda_u) = \lambda_{u1} + (\lambda_{ui}\cos\theta-\lambda_{uj}\sin\theta)i +
(\lambda_{ui}\sin\theta+\lambda_{uj}\cos\theta)j + \lambda_{uk}k.
$$
Consequently,
$$
\alpha(\lambda_u) = \lambda_u \Leftrightarrow \begin{pmatrix}\cos\theta & -\sin\theta\\
\sin\theta & \cos\theta \end{pmatrix} \begin{pmatrix} \lambda_{ui}\\\lambda_{uj}\end{pmatrix} =
\begin{pmatrix} \lambda_{ui}\\\lambda_{uj}\end{pmatrix}.
$$
So,
$$
\alpha(\lambda_u) = \lambda_u\, \text{ if and only if }\,\lambda_{ui} = \lambda_{uj} = 0.
$$
Then we assume that this condition is satisfied. Again by definition
\begin{align*}
&\lambda_u 1 = \alpha^u(1)\lambda_u,\\
&\lambda_u k = \alpha^u(k)(\lambda_{u1}+\lambda_{uk}k),\\
&\lambda_u i = (\lambda_{u1}+\lambda_{uk}k)i = \lambda_{u1}i+\lambda_{uk}j,\\
& \alpha^u(i)\lambda_u = (\cos(u\theta)i+\sin(u\theta)j)(\lambda_{u1}+\lambda_{uk}k) \\
& \phantom{\alpha^u(i)\lambda_u} = (\cos(u\theta)\lambda_{u1}+\sin(u\theta)\lambda_{uk})i +
(\sin(u\theta)\lambda_{u1}-\cos(u\theta)\lambda_{uk})j,\\
&\lambda_u j = (\lambda_{u1}+\lambda_{uk}k)j = \lambda_{u1}j-\lambda_{uk}i,\\
& \alpha^u(j)\lambda_u = (-\sin(u\theta)i+\cos(u\theta)j)(\lambda_{u1}+\lambda_{uk}k) \\
& \phantom{\alpha^u(j)\lambda_u} = (-\sin(u\theta)\lambda_{u1}+\cos(u\theta)\lambda_{uk})i +
(\cos(u\theta)\lambda_{u1}+\sin(u\theta)\lambda_{uk})j.
\end{align*}
Hence, $\alpha^u(\lambda)\lambda_u = \lambda_u\lambda$ for all $\lambda\in \mathbb{H}$ if and only
if
\begin{equation*}
\begin{pmatrix}\cos(u\theta) & \sin(u\theta)\\-\sin(u\theta) & \cos(u\theta) \end{pmatrix}
\begin{pmatrix} \lambda_{u1}\\\lambda_{uk}\end{pmatrix}= \begin{pmatrix} \lambda_{u1}\\
-\lambda_{uk}\end{pmatrix}.
\end{equation*}
Write
$$
\lambda_u = \lambda_{u1} + \lambda_{uk} k = \rho e^{k\beta} = \rho(\cos\beta + k\sin\beta)
$$
with $\rho\ge 0$ and $0\le\beta<2\pi$. The above matrix equality is equivalent to $\beta
=-u\theta/2 \pmod{\pi}$. Thus the following assertions are equivalent
\begin{itemize}

\smallskip

\item{} $\alpha(\lambda_u)= \lambda_u$ and $\alpha^u(\lambda)\lambda_u = \lambda_u\lambda$,

\smallskip

\item{} $\lambda_u = \rho e^{-ku\theta/2}$ with $\rho\in\mathbb{R}$.

\end{itemize}

\smallskip

Next we compute the complex $C_S(A)$ introduced in item~(1) of Theorem~\ref{teorema 3.1}.

\subsection{Computation of $A^{\alpha^r}$} Let $\gamma_0 + \gamma_1 x +\cdots +
\gamma_{n-1} x^{n-1}\in A$. Since
$$
\gamma_ux^u\lambda = \gamma_u \alpha^u(\lambda) x^u,
$$
we have that
$$
\gamma_0 + \gamma_1 x +\cdots + \gamma_{n-1} x^{n-1}\in A^{\alpha^r} \Leftrightarrow \gamma_u
\lambda = \alpha^{r-u}(\lambda)\gamma_u
$$
for all $u$ and all $\lambda \in \mathbb{H}$. Write $\gamma_u = \gamma_{u1} + \gamma_{ui}i +
\gamma_{uj}j + \gamma_{uk}k$. Since,
$$
\gamma_u k = \gamma_{u1}k - \gamma_{ui}j + \gamma_{uj}i - \gamma_{uk}\quad\text{and}\quad
k\gamma_u = \gamma_{u1}k + \gamma_{ui}j - \gamma_{uj}i - \gamma_{uk},
$$
if $\gamma_uk = \alpha^{r-u}(k)\gamma_u = k\gamma_u$, then $\gamma_u = \gamma_{u1} +
\gamma_{uk}k$. Arguing now for $\gamma_u$ as above for $\lambda_u$, we obtain that
$$
\gamma_u\lambda = \alpha^{r-u}(\lambda)\gamma_u\quad\text{for all $\lambda \in \mathbb{H}$,}
$$
if and only if
$$
\gamma_{u1}+\gamma_{uk}k = \rho e^{k(u-r)\theta/2} =
\rho\left(\cos\left(\frac{(u-r)\theta}{2}\right) + k\sin \left(\frac{(u-r)\theta}{2}\right)\right)
$$
with $\rho\in \mathbb{R}$. So, the following conditions are equivalent
\begin{itemize}

\smallskip

\item $\gamma_0 + \gamma_1 x +\cdots + \gamma_{n-1} x^{n-1}\in A^{\alpha^r}$.

\smallskip

\item $\gamma_u = \rho e^{k(u-r)\theta/2}$ with $\rho\in \mathbb{R}$.
\end{itemize}

\subsection{Computation of the boundary maps} By the above computations we know that if
$\gamma_n + \gamma_{n-1}x +\cdots + \gamma_1 x^{n-1}\in A^{\alpha^{mn}}$, then $\alpha(\gamma_u) =
\gamma_u$. Hence,
$$
d^{2m+1}(\gamma_n +\cdots + \gamma_1 x^{n-1}) = (\alpha(\gamma_n)-\gamma_n)x + \cdots +
(\alpha(\gamma_1)-\gamma_1) x^n = 0
$$
and
$$
d^{2m+2}(\gamma_n +\cdots + \gamma_1 x^{n-1}) = \left(\sum_{i=1}^n i\lambda_{n-i}
x^{i-1}\right)(\gamma_n +\cdots + \gamma_1 x^{n-1}).
$$

\subsection{Computation of the cohomology} Let $C = \mathbb{R}[x]/\langle g \rangle$, where $g =
x^n + \varsigma_1x^{n-1} +\cdots+\varsigma_n\in \mathbb{R}[x]$ is the polynomial defined by
$\lambda_u = \varsigma_u e^{-k\frac{u\theta}{2}}$. A direct computation, using the results
obtained in Subsections~6.1 and 6.2, shows that the map $\vartheta^*\colon C_S(C)\to C_S(A)$,
defined by
$$
\vartheta^{2m}(x^u) = e^{k(u-mn)\theta/2} x^u\qquad\text{and}\qquad \vartheta^{2m+1}(x^u) =
e^{k(u-mn-1)\theta/2} x^u,
$$
is an isomorphism of complexes. Since $\mathbb{H}$ is an $\mathbb{R}$-algebra separable
\begin{align*}
& \HH^0(A) = \HH^0(C) = C,\\
& \HH^{2m+1}(A) = \HH^{2m+1}(C) = \{h\in C:g'h = 0\},\\
& \HH^{2m+2}(A) = \HH^{2m+2}(C) = \frac{C}{\langle g' \rangle}.
\end{align*}
Moreover, from the formulas obtained in Theorem~3.2, it follows that $\vartheta^*$ induces an
isomorphism of algebras from $\HH_{\mathbb{R}}^*(C)$ to $\HH_{\mathbb{R}}^*(A)$. So, the product
of two homogeneous elements is induced by the multiplication map in $A$ whenever at least one of
them have even degree, and it is zero otherwise.

\end{document}